\documentclass{amsart}
\usepackage{amssymb,amsfonts,latexsym}

\newtheorem{theorem}{Theorem}[section]
\newtheorem{lemma}[theorem]{Lemma}
\newtheorem{proposition}[theorem]{Proposition}
\newtheorem{corollary}[theorem]{Corollary}
\theoremstyle{definition}
\newtheorem{definition}[theorem]{Definition}

\newtheorem{example}[theorem]{Example}

\newtheorem{remark}[theorem]{Remark}

\newcommand{\id}{\text{id}}

\newcommand{\Rep}{\text{Rep}}

\newcommand{\ot}{\otimes}

\newcommand{\ben}{\begin{enumerate}}
\newcommand{\een}{\end{enumerate}}

\begin{document}

\title[Frobenius-Schur indicators for semisimple Lie
algebras] {Frobenius-Schur indicators for semisimple Lie algebras}
\author{Mohammad Abu-Hamed}
\address{Department of Mathematics, Technion-Israel Institute of
Technology, Haifa 32000, Israel}
\email{mohammad@tx.technion.ac.il,
mohammad.abu-hamed@weizmann.ac.il}

\author{Shlomo Gelaki}
\address{Department of Mathematics, Technion-Israel Institute of
Technology, Haifa 32000, Israel}
\email{gelaki@math.technion.ac.il}

\date{February 11, 2007}

\maketitle

\section{Introduction}

 Classically Frobenius-Schur indicators were defined for irreducible
 representations of finite groups over the field of complex
 numbers. The interest in doing so came from the second indicator which
 determines whether an
 irreducible representation is real, complex or quaternionic. Namely,
 a classical theorem of Frobenius and Schur asserts that an irreducible
 representation is real, complex or quaternionic if and only if
 its second indicator is $1$, $0$ or $-1$, respectively (see e.g. \cite{s}). However,
 no representation-theoretic interpretation of the higher indicators is known.

 Recently, Frobenius-Schur indicators of irreducible representations of complex
 semisimple finite dimensional (quasi-)Hopf algebras $H$ were defined by
 Linchenko and Montgomery \cite{lm} and Mason and Ng \cite{mn} (see also
 \cite{ksz}), generalizing the
 definition in the group case. The values of the $m$th indicator are
cyclotomic integers in $\mathbb{Q}_{m}$. Moreover, an analog of
the Frobenius-Schur theorem on the second indicator was proved,
and in general it has been shown that the indicators carry rich
information on $H$, as well as on its representation category (see
also \cite{ns2}).

 In fact, one can generalize the definition of Frobenius-Schur indicators to simple objects
 of any semisimple tensor categories which admit a pivotal structure ($=$ tensor isomorphism
$\id\to ^{**}$), thus showing in particular that the indicators
are categorical invariants (see e.g. \cite{fgsv}, \cite{ns1}).

  The category of finite dimensional representations of a finite dimensional
  complex semisimple Lie algebra
  is a pivotal semisimple tensor category, and hence one can define the Frobenius-Schur indicators
  of its simple objects. The second indicator was already defined
  and known to be nonzero if and only if the simple representation is
self-dual, and $1$ or $-1$ if and only if the representation is
orthogonal or symplectic, respectively. Furthermore, Tits gave an
explicit formula for it in representation-theoretic terms (see
Section 3).

  The purpose of this paper is to study Frobenius-Schur indicators (of all
  degrees) for semisimple Lie algebras.
More specifically to find a closed formula for the indicators in
representation-theoretic terms and deduce its asymptotical
behavior. In particular we obtain that the indicators take integer
values.

The organization of the paper is as follows.

  Section 2 is devoted to preliminaries. We recall some basic
definitions and facts from Lie theory which we need (e.g. the Weyl
integration formula). Next we define the $m$th Frobenius-Schur
indicator of the representation categories of finite dimensional
complex semisimple Lie algebras.

  In section 3 we recall the properties of the second indicator. For the benefit
  of the reader we also give a proof of Tits' theorem.

Section 4 is dedicated to the proof of our main results. In 4.1 we
prove the formula for the $m$th Frobenius-Schur indicator
$\nu_{m}$, $m\geq 2$, which is given by the following theorem.

   \begin{theorem}\label{main}
 Let $\frak{g}$ be a finite dimensional complex
semisimple Lie algebra. Let $V(\lambda)$ be an irreducible
representation of
 $\mathfrak{g}$ with highest weight $\lambda$, $\mathcal{W}$ the Weyl group of $\frak{g}$, $\rho$
 the half sum of positive roots, and $V(\lambda)\left[\frac{\rho-\sigma\cdot\rho}{m}\right]$
 the weight space of the weight $\frac{\rho-\sigma\cdot\rho}{m}$ where $m\geq 2$ is an integer.
 Then the $m$th Frobenius-Schur indicator $\nu_{m}(V(\lambda))$ of
 $V(\lambda)$ is given by
    $$
    \displaystyle
    \nu_{m}(V(\lambda))=\sum_{\sigma\in\mathcal{W}}sn(\sigma)\
    dimV(\lambda)\left[\frac{\rho-\sigma\cdot\rho}{m}\right].
   $$

    \end{theorem}
Our proof of Theorem \ref{main} is analytic. Namely, we work with
the equivalent representation category of the associated simply
connected Lie group and use the Weyl integration formula to obtain
our formula.

Next, in 4.2 we prove the following corollary of Theorem
\ref{main}.

\begin{corollary}\label{conclusion}
For large enough $m$, $\nu_{m}(V(\lambda))=dimV(\lambda)[0]$
(which is not zero if and only if $\lambda$ belongs to the root
lattice). In particular for the classical Lie algebras
$\mathfrak{sl}(n,\mathbb{C})$, $\mathfrak{so}(2n,\mathbb{C})$,
$\mathfrak{so}(2n+1,\mathbb{C})$ and
$\mathfrak{sp}(2n,\mathbb{C})$,
$\nu_{m}(V(\lambda))=dimV(\lambda)[0]$ for $m$ greater or equal to
$2n-1$, $4n-5$, $4n-3$ and $2n+1$, respectively.
\end{corollary}

 Finally in 4.3 we use our formula and Kostant's theorem to compute explicitly
 the Frobenius-Schur indicators for the representation category of $sl(3,\mathbb{C})$.
 More specifically, we prove:

\begin{theorem}\label{sl3}
  Let $V(a,b)$ be an irreducible representation of $sl(3,\mathbb{C})$.
  Then
\begin{enumerate}
\item $\nu_{2}(V(a,b))=1$ if $a=b$, and $\nu_{2}(V(a,b))=0$ if $a\neq
b$.
\item $\nu_{3}(V(a,b))=1+min\{a,b\}$.
\item For $m>3$ we have, $\nu_{m}(V(a,b))=1+min\{a,b\}$ if
  $(a,b)$ is in the root lattice and $\nu_{m}(V(a,b))=0$ otherwise.
  \end{enumerate}
 \end{theorem}

\noindent{\bf Acknowledgments.} This research was supported by the
Israel Science Foundation (grant No. 125/05).

 \section{Preliminaries}

 Throughout let $\frak{g}$ be a finite dimensional complex semisimple Lie algebra
 of rank $r$, $(\,,\,)$ its Killing form, $\frak{h}$ a Cartan
 subalgebra (CSA) of $\frak{g}$, $\Phi$ the root system
 corresponding to $\frak{h}$, $\Delta$ a fixed base, $\{h_{1},...,h_{r}\}$ the
 corresponding coroot system, and $\mathcal{W}$ the Weyl group.

Let $\lambda\in\frak{h}^{*}$ be a dominant integral weight (i.e.
$\lambda(h_i)$ is a nonnegative integer for all $i$),
  $V(\lambda)$ the finite dimensional irreducible representation of $\frak{g}$ with
  highest weight $\lambda$ and $\Pi(\lambda)$ the set of integral weights occurring in
  $V(\lambda)$; it is a finite set which is invariant under the action of the Weyl group.
   For $\mu\in\Pi(\lambda)$, let $m_{\lambda}(\mu)=dimV(\lambda)[\mu]$ be the multiplicity
   of $\mu$ in $V(\lambda)$.
Recall that the multiplicities are invariant under the Weyl group
action.
 Let $ \displaystyle \rho=\frac{1}{2}\sum_{\alpha\in\Phi^{+}}\alpha$ (half sum of positive
 roots);
     it is a strongly dominant integral weight.

   Let us recall Kostant's theorem on the multiplicities of weights (for a proof see \cite{hu}).
    Let $\mu\in\frak{h}^{*}$ and define $p(\mu)$ to be
   the number of sets
   of non-negative integers $\{k_{\alpha}|\alpha\succ 0\}$ for
   which $\displaystyle \mu=\sum_{\alpha\succ 0}k_{\alpha}\alpha$ ($p$
   is called the Kostant's partition function). Of course,
   $p(\mu)=0$ if $\mu$ is not in the root lattice.

   \begin{theorem}\label{kostant}
   {\em (Kostant)} Let $\lambda$ be a dominant weight and $\mu\in\Pi(\lambda)$.
   Then the multiplicities of
  $V(\lambda)$ are given by the formula
  $$
 \displaystyle
  m_{\lambda}(\mu)=\sum_{\sigma\in\mathcal{W}}sn(\sigma)p(\sigma(\lambda+\rho)-\mu-\rho).
  $$
  \end{theorem}

Let $\frak{g}_{c}$ be the compact real form of $\frak{g}$, and $G$
the corresponding simply connected compact matrix Lie group with
Lie algebra $\frak{g}_{c}$. It is known that $\Rep(\frak{g})$,
$\Rep(\frak{g}_{c})$ and $\Rep(G)$ are equivalent symmetric tensor
categories.

 Let $\frak{t}$ be a CSA of $\frak{g}_{c}$; it corresponds to
 a maximal torus $T$ of $G$. Then
 $\frak{h}=\frak{t}\oplus i\frak{t}$.
 It is known that $\alpha(h)$ is purely imaginary for all $h\in\frak{t}$ and $\alpha\in\Phi$.
If $\frak{t}^{*}$ denotes the space of real-valued
       linear functionals on $\frak{t}$, then the roots are contained
       in $i\frak{t}^{*}\subset\frak{h}^{*}$.
       It is then convenient to introduce the {\em real} roots, which are simply
        $\frac{1}{i}$ times the ordinary roots, the {\em real} coroots $h_{\alpha}$
        which are the elements of $\frak{t}$ corresponding to the
        elements $\frac{2\alpha}{(\alpha,\alpha)}$ where
         $\alpha$ is a real root, and
        the {\em real}
        weights of an irreducible representation of $G$. An element $\mu$ of $\frak{t}^{*}$
        is said to be {\em integral} if $\mu(h_{\alpha})\in\mathbb{Z}$ for each
        real coroot $h_{\alpha}$. The real weights of any finite dimensional representation
        of $\frak{g}$ are integral. (See \cite{ha}.)

   The Weyl denominator is the
   function $A_{\rho}:T\longrightarrow\mathbb{C}$ given
   by $$
   A_{\rho}(t)=A_{\rho}(e^{h})=\sum_{\omega\in\mathcal{W}}sn(\omega)e^{i(\omega\cdot\rho)(h)}.
   $$

   \begin{theorem}{\em (Weyl integration formula)}
   Let $G$ be a simply connected compact Lie group. Let $f$ be a
   continuous class function on $G$,
   $dg$ the normalized Haar measure on $G$, and $dt$ the normalized
   Haar measure on $T$. Then
   $$
   \int_{G}f(g)dg=\frac{1}{|\mathcal{W}|}\int_{T}f(t)|A_{\rho}(t)|^{2}dt.
   $$
   \end{theorem}

   Let us now define the Frobenius-Schur indicators of an
   irreducible representation of $\frak{g}$.

 \begin{definition} Let $V$ be an irreducible representation of $\frak{g}$
   and $m\geq2$ be an integer. The $m$th Frobenius-Schur indicator of $V$ is the number
   $\nu_{m}(V)=tr(c|_{(V^{\otimes m})^{\frak{g}}})$, where
   $c$ is the cyclic automorphism of $V^{\otimes m}$ given by $v_1\ot\cdots\ot
   v_m\mapsto v_m\ot\cdots\ot v_{m-1}$.
   \end{definition}

\begin{remark}
   In fact, as we mentioned in the introduction, the indicators can be defined categorically.
   Applying the categorical definition to $\Rep(\frak{g})$ yields the above definition, while
   applying it to $\Rep(G)$ yields $tr(c|_{(V^{\otimes m})^{G}})$. Since the
   indicators of $V$ regarded as a $\frak{g}$-module coincide with the indicators of $V$
   regarded as a $G$-module we have
  $\nu_{m}(V)=tr(c|_{(V^{\otimes m})^{G}}).$
\end{remark}

\section{Tits' theorem on the second indicator}

\begin{theorem}\label{new} (See \cite{b})\label{f-s c}
Let $G$ be a compact Lie group.
  Let $V$ be an irreducible complex representation of $G$, and set $\epsilon_{V}=\int_{G}\chi(g^{2})dg$.
  Then $V$ is self dual if and only if $\epsilon_{V}\ne 0$.
  Furthermore, suppose $V$ is self dual and let $B$ be a (unique up to scalar) $G$-invariant
non-degenerate bilinear form on $V$. Then
  $B$ is either symmetric or skew-symmetric, and it is such if and
  only if $\epsilon_{V}=1,-1$, respectively.
  \end{theorem}

\begin{remark}
In Proposition \ref{f-s i} we will prove that
$\epsilon_{V}=\nu_{2}(V)$ as defined above. Historically
$\nu_{2}(V)$ was defined by $\epsilon_{V}$.
\end{remark}

\begin{example} Let us use Theorem \ref{main} to calculate $\nu_{m}(V)$ in the representation category
    of $sl(2,\mathbb{C})$.
    Let $sl(2,\mathbb{C})=sp\{h,x,y\}$,
    where $h= \left(
         \begin{matrix}
            1  & 0\\
            0 & -1
           \end{matrix}
           \right),$
$x= \left(
         \begin{matrix}
            0  & 0\\
            1 &  0
           \end{matrix}
           \right),$
$y= \left(
         \begin{matrix}
            0  & 1\\
            0 &  0
           \end{matrix}
           \right).$
 The root system is
    $\Phi=\{\alpha,-\alpha\}$, where
     $\alpha(h)=2$. The Weyl group is
    $\mathcal{W}=\{1,\sigma_{\alpha}\}$
     , $\rho=\frac{1}{2}\alpha$ and
     $\sigma_{\alpha}(\rho)=-\frac{1}{2}\alpha$.
     Let $V(n)=\oplus_{j=0}^{n}V[n-2j]$ be the irreducible representation of highest weight
     $\lambda(h)=n$ with its weight space decomposition.
     By Theorem \ref{main},
     $$
     \nu_{m}(V(n))=dimV(n)[0]-dimV(n)\left[\frac{\alpha}{m}\right].
     $$

  Let $m=2$. By the formula above, if $n$ is odd, then $dimV(n)[0]=0$ and
 $dimV(n)[\frac{\alpha}{2}]=1$. Hence $\nu_{2}(V(n))=-1$. Similarly, if $n$
  is even, $\nu_{2}(V(n))=1$.
 Consequently $\nu_{2}(V(n))=(-1)^{n}=(-1)^{\lambda(h)}$.

 For $m\geq 3$, $\frac{\alpha(h)}{m}=\frac{2}{m}$ is not an integer and hence
 $\nu_{m}(V(n))=dimV(n)[0]$. Therefore we have
     $$
     \nu_{m}(V(n))=\begin{cases}
                          \ \ 1 \ \ &\text{if n is even}\\
                          \ \ 0 \ \ &\text{if n is odd}
                         \end{cases}
    $$

 \end{example}

 Let $\frak{g}=\frak{h}\bigoplus (\displaystyle\oplus_{\alpha\in\Phi}\frak{g}_{\alpha})$
 be the root space decomposition of $\frak{g}$ and $\Delta =\{\alpha_{1},...,\alpha_{r}\}$ a fixed base.
 Fix a standard
 set of generators for $\frak{g}$: $x_{i}\in\frak{g}_{\alpha_{i}},\ \
 y_{i}\in\frak{g}_{-\alpha_{i}}$ so that $[x_{i},y_{i}]=h_{i}$.
Let $\check{\rho}:=1/2\sum_{\alpha\in \Phi^+} h_{\alpha}$ be the half sum of positive
 coroots.

 \begin{proposition} Let $E:=x_{1}+...+x_{r}$ and $H:=2\check{\rho}$. Then
 there exist constants $a_{1},...,a_{r}$ such that the subalgebra
 $P$ generated by $H,E,F:=a_{1}y_{1}+...+a_{r}y_{r}$ is isomorphic to $sl(2,\mathbb{C})$.
 \end{proposition}

 The Lie subalgebra $P\subseteq\frak{g}$ is called a principal $sl(2,\mathbb{C})$-subalgebra
 of $\frak{g}$ (see \cite{k} or \cite{dy}).

 \begin{lemma}\label{priciplesubalg}
 Let $V=V(\lambda)$ be an irreducible representation of $\frak g$.
  Let $P$ be a principal $sl(2,\mathbb{C})$-subalgebra of $\mathfrak{g}$.
   Consider $V$ as a P-module.
   Then its highest weight is $\lambda(H)$, and it
   contains the irreducible $sl(2,\mathbb{C})$-representation
   $V(\lambda(H))$ with multiplicity one.
  \end{lemma}

 \begin{proof}
 Let $v^{+}$ be a highest weight vector of $V$ considered
 as a $\mathfrak{g}$-module. Then obviously we have $Hv^{+}=\lambda(H)v^{+}$ and
 $Ev^{+}=0$. Hence $v^{+}$ is a highest weight vector with weight $\lambda(H)$
  for $V$ considered as a $P$-module.
  Therefore we can write $V=V(\lambda(H))\bigoplus V(n_{j})$. Now
 it remains to show that $\lambda(H)>n_{j}$ for any $j$.
 Let $\displaystyle V=V[\lambda]\oplus\bigoplus V[\mu]$
 be the weight space decomposition of $V$ as a $\frak{g}$-module. It is
 also a weight space decomposition of $V$ considered as a $P$-module,
  so $V[\mu]$ is a weight space of $P$ with weight $\mu(H)$. Recall
  that $\mu=\lambda-\sum_{j=1}^{r}k_{j}\alpha_{j}$ where
  $k_{j}\in\mathbb{Z}^{+}$. Note that
 $\lambda(H)>\mu(H)$ if and only if $\lambda(H)>
 \lambda(H)-\sum_{j=1}^{r}k_{j}\alpha_{j}(H)$ if and only if $
 \sum_{j=1}^{r}k_{j}\alpha_{j}(2\breve{\rho})>0$ if and only if $
 \sum_{j=1}^{r}k_{j}(\alpha_{j},2\rho)>0.$
 But $2\rho$ is strongly dominant,
 i.e., $(\alpha_{j},2\rho)>0$ for all $1\leq j \leq r$.
 The proof is complete.
 \end{proof}

Let $\omega_0\in \mathcal{W}$ be the unique element sending
$\Delta$ to $-\Delta$.

 \begin{theorem}\textbf{(Tits)}\label{tits}
 Let $V=V(\lambda)$ be a finite dimensional irreducible representation
 of $\frak{g}$. If $\lambda+\omega_{0}\lambda\neq
 0$ then $\nu_{2}(V)=0$. Otherwise,
 $\nu_{2}(V)=(-1)^{\lambda(2\check{\rho})}$.
 \end{theorem}

 \begin{proof}
 It is known that the dual of $V(\lambda)$ is $V(-\omega_{0}\lambda)$,
 so if
 $V(\lambda)$ is not self dual (i.e., $\lambda+\omega_{0}\lambda\neq 0$) then $\nu_{2}(V)=0$.

 Suppose that $V$ is self dual as a $\frak{g}$-module. Then $V$ admits a non-degenerate
 $\frak{g}$-invariant bilinear form, and we have to decide if it
 is symmetric or skew symmetric. To do so, consider the
 principal $sl(2,\mathbb{C})$-subalgebra $P$ as in Lemma \ref{priciplesubalg}. The
 restriction of $V$ to $P$ has a unique copy of the largest representation of $P$ occurring in
 $V$, with highest weight
 $\lambda(2\check{\rho})$. We already proved
 that this representation has indicator
 $(-1)^{\lambda(2\breve{\rho})}$. Now we can use Theorem \ref{new} to prove that $V$ has a symmetric
 (skew-symmetric) $\frak{g}$-invariant form if and only if it has a symmetric
 (skew-symmetric) $P$-invariant form. The first direction is
 obvious. Conversely, suppose that $V$ has a symmetric $P$-invariant
 form and suppose on the contrary that $V$ admits a skew-symmetric
 $\frak{g}$-invariant form. Then if we
 restrict the bilinear $\frak{g}$-form to $P$ we get that $V$ has a skew-symmetric $P$-invariant form
 which is a contradiction. Similar considerations are applied when $V$ has a skew-symmetric
 $P$-invariant form. We conclude that $\nu_{2}(V)=(-1)^{\lambda(2\breve{\rho})}$.
 \end{proof}

\section{The Main results}

 \subsection{Proof of Theorem \ref{main}}

   Let $G$ be the associated simply connected compact Lie group.
   \textit{ From now on we will consider $V(\lambda)$ as a $G$-module.}
For convenience set
  $V=V(\lambda)$, $N=V(\lambda)^{\otimes m}$, and let
  $\pi:G\longrightarrow GL(V)$ be the irreducible
  representation.

     The following lemma is easily derived from linear
     algebra.

       \begin{lemma}\label{lemma4.1} Let $T\in End(V)$ be a projection,
       $W=ImT$ and $S\in End(V)$ an operator preserving W.
       Then $tr|_{W}(S)=tr|_{V}(S\circ T).$
       \end{lemma}

       \begin{proof}
       Fix a basis $A=\{w_{1},...,w_{k}\}$ for $W$, and let
       $\tilde{A}=\{w_{1},...,w_{k},w_{k+1},...,w_{n}\}$ be a
       completion to a basis for $V$. Let $C=[S|_{W}]_{A}$ be the matrix representing
       $S|_{W}$ with respect to the basis $A$. Since $T|_{W}=id_{W}$ and $S(W)\subseteq W$ we
       find out that
        $
        [T]_{\tilde{A}}=
        \left(
         \begin{matrix}
           I_{k} & 0\\
             0   & 0
           \end{matrix}
           \right),
         $
        $
        [S]_{\tilde{A}}=
        \left(
         \begin{matrix}
           C   &  *\\
           0   &  *
           \end{matrix}
           \right),
         $
    and hence
    $
        [S]_{\tilde{A}}[T]_{\tilde{A}}=
        \left(
         \begin{matrix}
            C  & 0\\
            0  & 0
           \end{matrix}
           \right).
    $
  The lemma follows easily now.
  \end{proof}

  \begin{proposition}\label{pro 4.1}
  We have,
    $$
    \nu_{m}(V)=tr|_{N^{G}}(c)=\int_{G}tr|_{V}(c\circ\pi^{\otimes
   m}(g))dg.
    $$
   \end{proposition}

   \begin{proof}
   We follow the lines
   of the proof of the first formula for Frobenius-Schur
   indicators in the Hopf case, given
   in Section 2.3 of \cite{ksz}.

    Set $\tau=\pi^{\otimes m}$. Consider the operator $\int_{G}\tau(g)dg:N\longrightarrow
    N$.  Let us first
    show that the image of this operator is $N^{G}$. Indeed, by the invariance of the Haar
    measure,
    $\tau(h)\int_{G}\tau(g)vdg=\int_{G}\tau(hg)vdg=\int_{G}\tau(g)vdg$
    for all $h\in G$ and $v\in N$.
    Hence
    $
    Im\left(\int_{G}\tau(g)dg\right)\subseteq N^{G}.
    $

    Conversely, suppose that $u\in N^{G}$, then
    $
    \int_{G}\tau(g)udg=\int_{G}udg=u\int_{G}dg=u.
    $
     Hence $N^{G}\subseteq Im\left(\int_{G}\tau(g)dg\right)$ and we are done.

       In fact, the above shows also that the operator $\int_{G}\tau(g)dg$ is a projection
       onto $N^G$.

    Finally, $c\in Aut(N^{G})$, so by Lemma \ref{lemma4.1},
   $$tr|_{N^{G}}(c)=tr|_{N}\left(c\circ \int_{G}\tau(g)dg
    \right)=\int_{G}tr|_{N}(c\circ\tau(g))dg,
    $$
    as claimed.
 \end{proof}

The following lemma is a particular case of a lemma in Section 2.3
of \cite{ksz} and its proof replicates the proof of that lemma.

 \begin{lemma}\label{linalg}
  Let $f_{1},...,f_{m}\in
   End(V)$. Then,
   $$
   tr|_{V^{\otimes m}}(c\circ(f_{1}\otimes ...\otimes
   f_{m}))=tr|_{V}(f_{1}\circ ... \circ f_{m}) .
   $$
 \end{lemma}

 \begin{proof}
   Let $v_{1},...,v_{n}$ be a basis of V with dual basis $v^{*}_{1},...,v^{*}_{n}$.
   For $l=1,...,m$, $f_{l}$ is presented by the
   matrix $\left(a_{ij}^{l}\right)_{i,j=1}^{n}$, where
      $a_{ij}^{l}=(v_{i}^{*},f_{l}(v_{j}))$.
      Therefore,
   $tr\left(a_{ij}^{l}\right)=\sum_{i=1}^{n}(v_{i}^{*},f_{l}(v_{i}))$.
   We now have

$$\begin{array}{l}
  \displaystyle
    tr|_{V^{\otimes m}}(c\circ(f_{1}\otimes ... \otimes
   f_{m})) \ = \\
   \\
  \displaystyle\sum_{i_{1},...,i_{m}=1}^{n}(v_{i_{1}}^{*}\otimes v_{i_{2}}^{*}\otimes...\otimes
   v_{i_{m}}^{*}, c(f_{1}(v_{i_{1}})\otimes
   f_{2}(v_{i_{2}})\otimes...\otimes f_{m}(v_{i_{m}}))) \ =\\
\\
  \displaystyle\sum_{i_{1},...,i_{m}=1}^{n}(v_{i_{1}}^{*},f_{2}(v_{i_{2}}))\cdots
  (v_{i_{m-1}}^{*},f_{m}(v_{i_{m}}))
    (v_{i_{m}}^{*},f_{1}(v_{i_{1}})) \ =\\
\\
   \displaystyle\sum_{i_{1},...,i_{m}=1}^{n}
    a_{i_{1},i_{2}}^{2} a_{i_{2},i_{3}}^{3}\cdots a_{i_{m-1},i_{m}}^{m}
    a_{i_{m},i_{1}}^{1}=tr|_{V}(f_{2}\circ f_{3} \circ\cdots\circ f_{m}\circ f_{1}) \ = \\
\\
tr|_{V}(f_{1}\circ f_{2} \circ\cdots\circ f_{m}),
  \end{array}$$

   as desired.
\end{proof}

   Consequently we have the following proposition which is
   analogous to the finite group case.

  \begin{proposition}\label{f-s i}  Let $\chi$ be the irreducible
    character of $V$. Then
    $$
    \nu_{m}(V)=\int_{G}\chi(g^{m})dg.
    $$
  \end{proposition}

  \begin{proof}
We follow the lines
   of the proof of the first formula for Frobenius-Schur
   indicators in the Hopf case, given
   in Section 2.3 of \cite{ksz}.

    It follows immediately from Proposition \ref{pro 4.1} and Lemma \ref{linalg} that

    $$ \begin{array}{ll}
      \nu_{m}(V) = \int_{G}tr|_{N}(c\circ\pi^{\otimes
    m}(g))dg=\int_{G}tr|_{N}(c\circ(\pi(g)\otimes ...\otimes\pi(g))dg=
    \\ \\
        = \int_{G}tr|_{V}(\pi(g)\circ ..
    \circ\pi(g))dg=\int_{G}\chi(g^{m})dg. \\
    \end{array}$$
 \end{proof}

  Recall the integral real elements which are those elements $\mu$ of $\frak{t}^{*}$
  for which $\frac{2(\mu,\alpha)}{(\alpha,\alpha)}$ is an integer for
  any simple real root $\alpha$. For each real integral element $\mu$, there is a function
  $\tilde{\mu}$ on $T$ given by
  $$
  \tilde{\mu}(e^{h})=e^{i\mu(h)}
  $$
  for all $h$ in $\frak{t}$. Functions of this form are called
  torus characters and they have the following property.

\begin{lemma}\label{torusc}
$$
 \int_{T}\tilde{\mu}(t)dt=\int_{T}e^{i\mu(h)}de^{h}=
  \begin{cases}
    1 & \text{$\mu=0$}, \\
    0 & \text{otherwise}.
  \end{cases}
$$
\end{lemma}

\begin{proof}
 Suppose that $\mu\neq 0$, then there exists $t_{0}\in\frak{t}$ such
 that $\tilde{\mu}(t_{0})\neq 1$. Therefore
 $$
 \int_{T}\tilde{\mu}(t)dt=\int_{T}\tilde{\mu}(t_{0}t)dt
 =\tilde{\mu}(t_{0})\int_{T}\tilde{\mu}(t)dt,
 $$
 hence $\int_{T}\tilde{\mu}(t)dt=0$. \end{proof}

  Let $\chi$ be the character of $V$. Before we begin the proof of
  Theorem \ref{main}, recall that if $t=e^{h}\in T$ then for all $t\in T$,
    \begin{equation}
    \chi(t)=\chi(e^{h})=\sum_{\mu\in\Pi(V)}dim(V[\mu])e^{i\mu(h)}.
    \end{equation}

 We can now prove our main result.\\

  \noindent \textbf{Proof of Theorem \ref{main}}:
    By Proposition \ref{f-s i} and the Weyl integration formula we have,
    \begin{equation}\label{eq 2}
              \nu_{m}(V)=\int_{G}\chi(g^{m})dg=\frac{1}{|\mathcal{W}|}\int_{T}\chi(t^{m})|A_{\rho}(t)|^{2}dt.\\
      \end{equation}
      On the other hand,
    \begin{equation}\label{eq 3}
      \chi(t^{m})=\chi(e^{mh})=\sum_{\mu\in\Pi(V)}dim(V[\mu])e^{im\mu(h)}. \\
    \end{equation}
    Hence by (\ref{eq 2}) and (\ref{eq 3}) we have,
    \begin{equation}\label{eq 4}
    \nu_{m}(V)=\frac{1}{|\mathcal{W}|}\sum_{\mu\in\Pi(V)}dimV[\mu]
    \int_{T}e^{im\mu(h)}|A_{\rho}(e^{h})|^{2}de^{h}.\\
    \end{equation}

    Now let us calculate the last integral.
We have
    $$\begin{array}{ll}
    \int_{T}e^{im\mu(h)}|A_{\rho}(e^{h})|^{2}de^{h}=
    \int_{T}e^{im\mu(h)}A_{\rho}(e^{h})\overline{A_{\rho}(e^{h})}de^{h}=\\
    \\
    =\int_{T}e^{im\mu(h)}\left( \sum_{\omega\in
    \mathcal{W}}sn(\omega)e^{i(\omega\cdot\rho)(h)}\right)\left( \sum_{\tau\in
    \mathcal{W}}sn(\tau)e^{-i(\tau\cdot\rho)(h)}\right)de^{h}=\\ \\
    =\sum_{\omega,\tau\in
    \mathcal{W}}sn(\omega\tau)\int_{T}e^{i(m\mu+\omega \cdot\rho-\tau\cdot\rho)(h)}de^{h}. \\
    \end{array}$$

    But from Lemma \ref{torusc} we have
    \begin{equation}
     \int_{T}e^{i(m\mu+\omega
     \cdot\rho-\tau\cdot\rho)(h)}de^{h}=
       \begin{cases}
         1     &\text{if $m\mu+\omega\cdot\rho-\tau\cdot\rho= 0$}\\

         0     &\text{otherwise.}
       \end{cases}
    \end{equation}
       Hence (\ref{eq 4}) becomes,
    \begin{eqnarray*}
     \lefteqn {\nu_{m}(V)=\frac{1}{|\mathcal{W}|}\sum_{\omega,\tau\in\mathcal{W}}
     \sum_{\mu=\frac{\tau\cdot\rho-\omega\cdot\rho}{m}}
     sn(\omega\tau)dimV[\mu]}\\
     & = & \frac{1}{|\mathcal{W}|}\sum_{\omega
     ,\tau\in\mathcal{W}}
     sn(\omega\tau)dimV\left[\frac{\tau\cdot\rho-\omega\cdot\rho}{m}\right].
     \end{eqnarray*}
     Since $dimV[\zeta]=dimV[\tau\cdot\zeta]$ for all $\zeta\in\Pi(V)$ and
     $\tau\in\mathcal{W}$,
      we can write,
      $$
      \nu_{m}(V)=
      \frac{1}{|\mathcal{W}|}\sum_{\omega ,\tau\in\mathcal{W}}sn(\omega\tau)
      dimV\left[\frac{\rho-\tau^{-1}\omega\cdot\rho}{m}\right].
      $$

      Now if we fix $\omega\in\mathcal{W}$, substitute $\sigma=\tau^{-1}\omega$
      and use the fact that $sn(\omega\tau)=sn(\tau^{-1}\omega)$, we get
      \begin{gather*}
        \sum_{\tau\in\mathcal{W}}sn(\omega\tau)dimV\left[\frac{\rho-\tau^{-1}\omega\cdot\rho}{m}\right]=
   \sum_{\tau\in\mathcal{W}}sn(\tau^{-1}\omega)dimV\left[\frac{\rho-\tau^{-1}\omega\cdot\rho}{m}\right]=\\
    \sum_{\sigma\in\mathcal{W}}sn(\sigma)dimV\left[\frac{\rho-\sigma\cdot\rho}{m}\right].
      \end{gather*}

       Consequently,
      \begin{gather*}
     \nu_{m}(V)=\frac{1}{|\mathcal{W}|}\sum_{\omega,\sigma\in\mathcal{W}}
     sn(\sigma)dimV\left[\frac{\rho-\sigma\cdot\rho}{m}\right]
     =\sum_{\sigma\in\mathcal{W}}sn(\sigma)dimV\left[\frac{\rho-\sigma\cdot\rho}{m}\right],
     \end{gather*}
      as desired. \qed

It may be interesting to state the following immediate consequence
of Theorem \ref{main} and Theorem \ref{tits}.

 \begin{corollary}
  Let $V(\lambda)$ be an irreducible self dual representation of $\frak{g}$,
  then
  $$
 \displaystyle \sum_{\sigma\in\mathcal{W}}sn(\sigma)dimV(\lambda)\left[\frac{\rho-\sigma\cdot\rho}{2}\right]
 =(-1)^{\lambda(2\check{\rho})}.
  $$
If $V(\lambda)$ is not self dual, the sum equals $0$.
 \end{corollary}

 \subsection{Proof of Corollary \ref{conclusion}}

     Since $\rho$ is strongly dominant, $\sigma\cdot\rho=\rho$ only when
     $\sigma=1$. Write
     $$
      \nu_{m}(V)=dimV[0]+\sum_{\sigma\neq1}sn(\sigma)dimV
      \left[\frac{\rho-\sigma\cdot\rho}{m}\right]
     .$$
    We wish to show that for large enough $m$,
    $\frac{\rho-\sigma\cdot\rho}{m}$ is not a weight of $V$
    when $\sigma\neq 1$. Indeed, suppose that $\sigma\neq 1$.
     Recall that $\rho-\sigma\cdot\rho$ is an integral element, hence if we fix some coroot
     $h_{\alpha}$, we have the following
     set of integers: $U_{\alpha}=\{(\rho-\sigma\cdot\rho)(h_{\alpha})|
     \sigma\in\mathcal{W}, \sigma\neq 1\}$. Therefore if we take
     $m_{\alpha}=1+u_{\alpha}$, where $u_{\alpha}$ is the maximal
     element of $U_{\alpha}$, then
     $\frac{\rho-\sigma\cdot\rho}{m_{\alpha}}\notin \Pi(V)$.
      Hence $dimV\left[\frac{\rho-\sigma\cdot\rho}{m_{\alpha}}\right]=0$ for
     all $\sigma\neq 1$, and therefore
     $\nu_{m}(V)=dimV[0]$, for all $m\geq m_{\alpha}$. \qed

     Note that by the procedure of the above proof,
     $m=:min\{m_{\alpha}| \alpha\in\Delta\}$ is a better
     bound. Let us now give an explicit such lower bound.

     \begin{lemma}
     If $\omega\in\mathcal{W}$ then
     $$
     \displaystyle
     \omega\cdot\rho=\rho-\sum_{\substack{\alpha\in\Phi^{+} \\ \omega^{-1}(\alpha)\in\Phi^{-}}} \alpha.
     $$
      In particular, $ s_{\alpha}(\rho)=\rho-\alpha$ for
      $\alpha\in\Delta$.
      \end{lemma}

      \begin{proof} Evidently, $\omega\cdot\rho$ is half sum of the set
       $\{\omega(\alpha)|\alpha\in\Phi^{+}\}$. Like $\Phi^{+}$,
       this is a set of exactly half of the roots, containing
       each root or its negative but not both. More precisely,
       this set is obtained from $\Phi^{+}$ by replacing each
       $\alpha\in\Phi^{+}$ such that
       $\omega^{-1}\cdot\alpha\in\Phi^{-}$ by its negative. Now,
     $$
     \displaystyle
     \omega\cdot\rho=\rho-\sum_{\substack{\alpha\in\Phi^{+} \\ \omega^{-1}(\alpha)\in\Phi^{-}}} \alpha
     $$
     is evident , and $s_{\alpha}(\rho)=\rho-\alpha$ is a
     special case since one shows that if $\alpha\in\Delta$
     and $\beta\in\Phi^{+}$, then either $\beta=\alpha$ or
     $s_{\alpha}(\beta)\in\Phi^{+}$.
     \end{proof}

     \begin{proposition}\label{bound} Let $V$ be an irreducible representation of
      $G$. Then $\nu_{m}(V)=dimV[0]$ for all
      $ m\geq M:=\displaystyle
      min_{\alpha\in\Delta}\{\sum_{\beta\in\Phi^{+}}|\beta(h_{\alpha})|+1\}$.
      \end{proposition}

 \begin{proof}
     Let $h=h_{\alpha}$ be a simple coroot. For all $1\ne \omega\in\mathcal{W}$ we have,
     \begin{equation*}
     \left|(\rho-\omega\cdot\rho)(h)\right|=\left|\sum_{\substack{\beta\in\Phi^{+}, \\
     \omega^{-1}(\beta)\in\Phi^{-}}}
     \beta(h)\right|\leq \sum_{\substack{\beta\in\Phi^{+}, \\ \omega^{-1}(\beta)\in\Phi^{-}}}
     |\beta(h)|\leq\sum_{\beta\in\Phi^{+}}|\beta(h)|.
    \end{equation*}
    Therefore if we choose $m=\sum_{\beta\in\Phi^{+}}|\beta(h)|+1$ then
   $\frac{\rho-\omega\cdot\rho}{m}(h) \notin \mathbb{Z}$,
   namely, $\frac{\rho-\omega\cdot\rho}{m}$ is not a weight.
   Consequently,
   $\displaystyle \sum_{\sigma\neq1}sn(\sigma)dimV
   \left[\frac{\rho-\sigma\cdot\rho}{m}\right]=0$, and we are done.
   \end{proof}

    Let us calculate the bound $M$ defined in Proposition \ref{bound}
    for $\mathfrak{sl}(n,\mathbb{C})$. Let
    the Cartan subalgebra be the set of diagonal
    matrices in $\mathfrak{sl}(n,\mathbb{C})$.
Let the set of positive roots be $\Phi^+=\{\beta_{i,j}| 1\leq
i<j\leq n\}$, where $\beta_{i,j}(diag(a_1,\dots,a_n))=a_i-a_j$.
The subset $\Delta=\{\beta_{i,i+1}| 1\leq i\leq n-1\}$ is a base.
With respect to this base the simple coroots are $\{h_i|1\leq
i\leq n-1\}$, where $h_i$ is the matrix with $1$ in the $(i,i)$
position, $-1$ in the $(i+1,i+1)$ position and $0$ elsewhere.
Then, by an elementary calculation, we get that for any simple
coroot $h$,
   \begin{equation*}
   \sum_{1\leq i<j\leq n}|\beta_{i,j}(h)|+1=2n-1.
   \end{equation*}
Consequently we obtain that $M=2n-1$.

Let us calculate the bound $M$ defined in Proposition \ref{bound}
for $\mathfrak{so}(2n+1,\mathbb{C})$. Let the Cartan subalgebra be
the set of diagonal matrices in $\mathfrak{so}(2n+1,\mathbb{C})$.
Let the set of positive roots be $\Phi^+=\{\beta_{i}\pm
\beta_{j}|1\leq i<j\leq n\}\cup \{\beta_{i}|1\leq i\leq n\}$,
where $\beta_i(h_j)=\delta_{ij}$. The subset
$\Delta=\{\beta_{i}-\beta_{i+1},\,\beta_n| 1\leq i\leq n-1\}$ is a
base. With respect to this base the simple coroots are $\{h_{i}-
h_{i+1}, 2h_{n}|1\leq i\leq n-1\}$, where $h_i$ is the matrix with
$1$ in the $(i,i)$ position, $-1$ in the $(n+i,n+i)$ position and
$0$ elsewhere. Then, by an elementary calculation, we get that for
any simple coroot $h:=h_{k}-h_{k+1}$, $1\leq k\leq n-1$, the sum
$\sum_{\beta\in\Phi^{+}}|\beta (h)|+1$ equals
\begin{equation*}
\sum_{1\leq i<j\leq n}|(\beta_{i}+\beta_{j})(h)|+\sum_{1\leq
i<j\leq n}
|(\beta_{i}-\beta_{j})(h)|+\sum_{i=1}^{n}|\beta_{i}(h)|+1=4n-3,
\end{equation*}
while for the simple coroot $h:=2h_{n}$ it equals $4n-1$.
Consequently we obtain that $M=4n-3$.

Applying similar arguments to the other classical simple Lie
algebras yields the following result.

\begin{proposition}
The bound $M$ for $\mathfrak{sl}(n,\mathbb{C})$,
$\mathfrak{so}(2n,\mathbb{C})$, $\mathfrak{so}(2n+1,\mathbb{C})$
and $\mathfrak{sp}(2n,\mathbb{C})$ is equal to $2n-1$, $4n-5$,
$4n-3$ and $2n+1$, respectively.
\end{proposition}

\subsection{The proof of Theorem \ref{sl3}}

 Let $\frak{h}$ be the CSA of $sl(3,\mathbb{C})$ generated by the two elements
 $h_{1}=diag(1,-1,0)$ and $h_{2}=diag(0,1,-1)$.
 We will identify any functional $\alpha$ on $\frak{h}$ with the pair
 $(\alpha(h_{1}),\alpha(h_{2}))$. Under this identification the six roots of
 $sl(3,\mathbb{C})$ are
 $\alpha_{1}=(2,-1)$,
  $\alpha_{2}=(-1,2)$,
  $\alpha_{1}+\alpha_{2}=(1,1)$,
  $-\alpha_{1}=(-2,1)$,
  $-\alpha_{2}=(1,-2)$ and
  $-\alpha_{1}-\alpha_{2}=(-1,-1)$.
 The roots $\alpha_{1}=(2,-1),\ \alpha_{2}=(-1,2)$ form
 a base and the corresponding simple
 coroots are
 $h_{1},h_{2}$, respectively.

 Recall that if $V=V(\lambda)$ is an irreducible
 representation of $sl(3,\mathbb{C})$ of highest weight $\lambda$,
 then $\lambda$ is of the form $(a,b)$ with $a$ and $b$
 non-negative integers.

 Recall that $\mathcal{W}\cong S_{3}$ and it acts on $\frak{h}$ by
$\sigma\cdot diag(d_1,d_2,d_3)=
           diag(d_{\sigma(1)},d_{\sigma(2)},d_{\sigma(3)})$.
 Therefore, $(12)\cdot \alpha_{1}=-\alpha_{1}$, $(12)\cdot \alpha_{2}=
\alpha_{1}+\alpha_{2}$; $(13)\cdot \alpha_{1}=-\alpha_{2}$,
$(13)\cdot \alpha_{2}=-\alpha_{1}$; $(23)\cdot
\alpha_{1}=\alpha_{1}+\alpha_{2}$, $(23)\cdot
\alpha_{2}=-\alpha_{2}$; $(123)\cdot
\alpha_{1}=-\alpha_{1}-\alpha_{2}$, $(123)\cdot
\alpha_{2}=\alpha_{1}$; and $(132)\cdot \alpha_{1}=\alpha_{2}$,
$(123)\cdot \alpha_{2}=-\alpha_{1}-\alpha_{2}$.

The half sum of positive roots is
 $\rho=\frac{1}{2}(2\alpha_{1}+2\alpha_{2})=\alpha_{1}+\alpha_{2}$.
 We have, $\rho-(1 2)\rho=(2,-1)$,
 $\rho-(1 3)\rho=(2,2)$,
 $\rho-(2 3)\rho=(-1,2)$,
 $\rho-(1 2 3)\rho=(0,3)$, and
 $\rho-(1 3 2)\rho=(3,0)$.

 Let $m=2$. Considering our formula, we cancel all the summands which include roots
 that one of their two components is not divisible by 2.
 Consequently we get $$\nu_{2}(V)=dimV[(0,0)]-dimV[(1,1)].$$
Recall that an irreducible representation $V(a,b)$ is self dual if
and only if
 $a=b$. Since
 $\lambda=(s,s)=s\alpha_{1}+s\alpha_{2}$,
 $(\lambda,2\check{\rho})=(\lambda,2h_{1}+2h_{2})=
 (\lambda,2h_{1})+(\lambda,2h_{2})=4s$, it follows from Tits'
 theorem that
  $$
  \nu_{2}(V(a,b))=\begin{cases}
    0 & \text{$a\neq b$}, \\
    1 & \text{otherwise}.
  \end{cases}
$$

 Similar considerations for $m\geq 3$ yield,
$$
 \nu_{3}(V)=dimV[(0,0)]+dimV[(1,0)]+dimV[(0,1)]
$$
 and
 $$
 \nu_{m\geq 4}(V)=dimV[(0,0)].
 $$
In particular, if $\lambda$ does not belong to the root lattice,
$\nu_{m\geq 4}(V)=0$.

We now calculate
 $dimV[(0,0)]$, $dimV[(1,0)]$ and $dimV[(0,1)]$.
 Recall that for $\eta\in \mathfrak{h}^*$, $p(\eta)\ge 1$ if and only if
 $\eta$ belongs to the root lattice and $\eta\succ 0$. If $\eta=k\alpha_{1}+l\alpha_{2}$
 with nonnegative integers $k$ and $l$,
 then $p(\eta)=1+min\{k,l\}$. Write $\lambda=k\alpha_{1}+l\alpha_{2}$
 where $k$ and $l$ are real numbers and
 identify it with the pair
 $(\lambda(h_{1}),\lambda(h_{2}))=(a,b)=(2k-l,2l-k)$.

Note that $(0,1)=\frac{1}{3}\alpha_{1}+\frac{2}{3}\alpha_{2}$ and
$(1,0)=\frac{2}{3}\alpha_{1}+\frac{1}{3}\alpha_{2}$. Therefore by
Kostant's formula (see Theorem \ref{kostant}),
$$
dimV[(0,0)]=
\sum_{\omega\in\mathcal{W}}sn(\omega)p((k+1)\omega\cdot\alpha_{1}+(l+1)\omega
\cdot\alpha_{2}-\alpha_{1}-\alpha_{2}),
$$
$$
dimV[(0,1)]=
\sum_{\omega\in\mathcal{W}}sn(\omega)p\left((k+1)\omega\cdot\alpha_{1}+(l+1)\omega\cdot\alpha_{2}-
\frac{4}{3}\alpha_{1}-\frac{5}{3}\alpha_{2}\right),
$$
and
$$
dimV[(1,0)]=
\sum_{\omega\in\mathcal{W}}sn(\omega)p\left((k+1)\omega\cdot\alpha_{1}+(l+1)\omega\cdot\alpha_{2}-
\frac{5}{3}\alpha_{1}-\frac{4}{3}\alpha_{2}\right).
$$

It is straightforward to verify that in each of the three cases
the surviving terms correspond to $\omega=1,(12),(23)$. For
example, in the first case calculating \linebreak
$(k+1)\omega\cdot\alpha_{1}+(l+1)\omega
\cdot\alpha_{2}-\alpha_{1}-\alpha_{2}$
  for $\omega=(12),(13),(23),(123),(132)$, yields $(l-k-1)\alpha_{1}+l\alpha_{2}$,
  $-(l+2)\alpha_{1}-(k+2)\alpha_{2}$ (hence $p=0$),
  $k\alpha_{1}+(k-l-1)\alpha_{2}$,
  $(l-k-1)\alpha_{1}-(k+2)\alpha_{2}$ (hence $p=0$), and
  $-(l+2)\alpha_{1}+(k-l-1)\alpha_{2}$ (hence $p=0$), respectively.

Therefore we have that $dimV[(0,0)]$ equals
$$
 p\left(\left(\frac{b+2a}{3}\right )\alpha_{1}+\left(\frac{2b+a}{3}\right)\alpha_{2}\right) -
 p\left(\left(\frac{b-a-3}{3}\right
 )\alpha_{1}+\left(\frac{2b+a}{3}\right)\alpha_{2}\right)
 $$
 $$
- p\left(\left(\frac{b+2a}{3}\right )\alpha_{1}+
 \left(\frac{-b+a-3}{3}\right )\alpha_{2}\right ),
$$
$dimV[(0,1)]$ equals
$$
p\left(\left(\frac{b+2a-1}{3}\right
)\alpha_{1}+\left(\frac{2b+a-2}{3}\right)\alpha_{2}\right)-
 p\left(\left(\frac{b-a-4}{3}\right
 )\alpha_{1}+\left(\frac{2b+a-2}{3}\right)\alpha_{2}\right)
 $$
 $$
 - p\left(\left(\frac{b+2a-1}{3}\right )\alpha_{1}+
 \left(\frac{-b+a-5}{3}\right )\alpha_{2}\right ),
$$
and $dimV[(1,0)]$ equals
$$
p\left(\left(\frac{b+2a-2}{3}\right
)\alpha_{1}+\left(\frac{2b+a-1}{3}\right)\alpha_{2}\right)-
 p\left(\left(\frac{b-a-5}{3}\right
 )\alpha_{1}+\left(\frac{2b+a-1}{3}\right)\alpha_{2}\right)
 $$
 $$
-  p\left(\left(\frac{b+2a-2}{3}\right )\alpha_{1}+
 \left(\frac{-b+a-4}{3}\right )\alpha_{2}\right ).
$$

Now, modulo 3, exactly one of the following holds: 1) $b+2a=0$ and
$2b+a=0$ (in this case $\lambda$ belongs to the root lattice), 2)
$b+2a=1$ and $2b+a=2$ and 3) $b+2a=2$ and $2b+a=1$. Hence by the
above and elementary calculations, we obtain that in the first
case $dimV[(0,1)]=dimV[(1,0)]=0$, in the second case
$dimV[(0,0)]=dimV[(1,0)]=0$ and in the third case
$dimV[(0,0)]=dimV[(0,1)]=0$. Therefore, in the first case
$\nu_{3}(V(a,b))$ equals
$$
max\left \{1+min\left
\{\frac{b+2a}{3},\frac{2b+a}{3}\right\},
  0\right \}-  max \left \{1+min\left
\{\frac{b-a-3}{3},\frac{2b+a}{3}\right \},0 \right \}
$$
$$
- max\left \{1+min\left \{\frac{b+2a}{3},\frac{-b+a-3}{3}\right
\},0 \right\},
$$
in the second case it equals
$$
max\left \{1+min\left
\{\frac{b+2a-1}{3},\frac{2b+a-2}{3}\right\},
  0\right \}
- max \left \{1+min\left \{\frac{b-a-4}{3},\frac{2b+a-2}{3}\right
\},0 \right \}
$$
$$
- max\left \{1+min\left \{\frac{b+2a-1}{3},\frac{-b+a-5}{3}\right
\},0 \right\},
$$
and in the third case it equals
$$
 max\left \{1+min\left \{\frac{b+2a-2}{3},\frac{2b+a-1}{3}\right\},
  0\right \}
  - max \left \{1+min\left
\{\frac{b-a-5}{3},\frac{2b+a-1}{3}\right \},0 \right \}
$$
$$
-max\left \{1+min\left \{\frac{b+2a-2}{3},\frac{-b+a-4}{3}\right
\},0 \right\}.
$$

Finally, it is easy to check that in each case the sum equals
$1+min\{a,b\}$, as claimed. This completes the proof of the
theorem. \qed


\begin{thebibliography}{ABCD}

\bibitem[B]{b}
D. Bump, \textit{Lie Groups}, Springer-Verlag NY, LLC, (2004).

\bibitem[D]{dy}
 E. Dynkin, Semisimple subalgebras of semisimple Lie algebras
 (Russian) Mat.Sbornik N.S. \textbf{30 (27)} (1952) 349-462,
 English: \textbf{AMS} Translations \textbf{6} (1957), 111-244.

\bibitem[FGSV]{fgsv}
J. Fucs, C. Ganchev, K. Szlach$\acute{a}$nyi, and P. Vescernyes,
\textit{$S_{4}$-symmetry of 6j-sympols and Frobenius-Schur
indicators in rigid monoidal $\mathcal{C}^{*}$-categories}. J.Math
Phys. 40 (1999), 408-426.

\bibitem[Ha]{ha} B. Hall, \textit{Lie groups, Lie
algebras and representations}, Springer-Verlag,
Berlin-Heidelberg-New York, (2006).

\bibitem[Hu]{hu} J. Humphreys, \textit{Introdution to Lie
algebras and representation theory}, Springer-Verlag,
Berlin-Heidelberg-New York, (1972).

\bibitem[K]{k}
B. Kostant, \textit{The principal three dimensional subgroup and
betti numbers of complex simple Lie group}, Amer.J.Math. 81
(1959), 973-1032.

\bibitem[KSZ]{ksz}
Y. Kashina, Y. Sommerhaeuser, and Y. Zhu, \textit{On higher
Frobenius-Schur indicators}, Memoirs of the AMS 181, no 855
(2006).

\bibitem[LM]{lm}
V. Linchenko and S. Montgomery, \textit{A Frobenius-Schur theorem
for Hopf algebras}, Algebr. Represent. Theory 3 (2000), no. 4,
347-355, Special issue dedicated to Klaus Roggenkamp on the
occasion of his 60th birthday.

\bibitem[MN]{mn} G. Mason and S-H. Ng, \textit{Central invariants
and Frobenius-Schur indicators for semisimple qusi-Hopf algebras},
Adv. Math. 190 (2005), 161-195.

\bibitem[NS1]{ns1}
S-H. Ng and P. Schauenburg, \textit{Higher Frobenius-Schur
indicators for pivotal categories}, preprint
\textbf{arXiv:math,QA/0503167}.

\bibitem[NS2]{ns2}
S-H. Ng and P. Schauenburg, \textit{Central invariants and higher
indicators for semisimple quasi-Hopf algebras}, Transactions of
the AMS, to appear, \textbf{arXiv:math,QA/0508140}.

\bibitem[S]{s} J-P. Serre, \textit{Linear Representation of
Finite Groups}, Springer-Verlag, New York, (1977).

\end{thebibliography}
 \end{document}